\documentclass[a4paper,12pt,reqno]{amsart}
\usepackage[T1]{fontenc}
\usepackage[utf8]{inputenc}
\usepackage[libertinus,vvarbb]{newtx}
\usepackage[margin=1in]{geometry}
\usepackage{enumitem}
\usepackage{graphicx}
\usepackage[svgnames]{xcolor}
\usepackage{xparse}
\usepackage{xurl}
\usepackage[bookmarksdepth=2]{hyperref}

\hypersetup{colorlinks=true,urlcolor=MidnightBlue,linkcolor=MidnightBlue,citecolor=MidnightBlue}

\newtheorem{theorem}{Theorem}

\theoremstyle{definition}

\theoremstyle{remark}

\DeclareMathOperator{\id}{Id}
\DeclareMathOperator{\re}{Re}

\newcommand{\C}{\mathbb{C}}
\newcommand{\R}{\mathbb{R}}
\newcommand{\sph}{\mathbb{S}}
\newcommand{\ball}{\mathbb{B}}
\newcommand{\schwartz}{\mathscr{S}}
\newcommand{\fourier}{\mathscr{F}}
\newcommand{\ph}{\varphi}

\renewcommand{\le}{\leqslant}
\renewcommand{\ge}{\geqslant}

\newcommand{\lv}{\lvert}
\newcommand{\rv}{\rvert}
\newcommand{\lV}{\lVert}
\newcommand{\rV}{\rVert}

\NewDocumentCommand{\formula}{ssom}{%
 \IfBooleanTF{#1}{%
  \IfBooleanTF{#2}{%
   \IfValueTF{#3}%
    {\begin{align}\label{#3}\begin{gathered}#4\end{gathered}\end{align}}%
    {\begin{gather}#4\end{gather}}%
  }{%
   \IfValueTF{#3}%
    {\begin{align}\label{#3}\begin{aligned}#4\end{aligned}\end{align}}%
    {\begin{gather*}#4\end{gather*}}%
  }%
 }{%
  \IfValueTF{#3}%
   {\begin{align}\label{#3}#4\end{align}}%
   {\begin{align*}#4\end{align*}}%
 }%
}

\newcommand{\nist}[2]{\href{https://dlmf.nist.gov/#1\#E#2}{Eq.~#1.#2}}
\newcommand{\doi}[1]{\href{https://doi.org/#1}{\textsf{\scriptsize DOI:#1}}}
\newcommand{\arxiv}[1]{\href{https://arxiv.org/abs/#1}{\textsf{\scriptsize arXiv:#1}}}
\newcommand{\isbn}[1]{\textsf{ISBN:#1}}
\NewDocumentCommand{\link}{oom}{\href{#3}{\textsf{\scriptsize \IfValueTF{#1}{#1}{#3}}}\IfValueT{#2}{{\textsf{\scriptsize #2}}}}

\begin{document}

\title[Maximal inequalities for derivatives of spherical means]{Maximal inequalities for derivatives of spherical means}
\author{Mateusz Kwaśnicki}
\thanks{Work supported by the National Science Centre, Poland, grant no.\@ 2023/49/B/ST1/04303}
\address{Mateusz Kwaśnicki \\ \textnormal{Department of Analysis and Stochastic Processes \\ Wrocław University of Science and Technology \\ Wybrzeże Wyspiańskiego 27 \\ 50-370 Wrocław, Poland}}
\email{\href{mailto:mateusz.kwasnicki@pwr.edu.pl}{\textsf{mateusz.kwasnicki@pwr.edu.pl}}}

\begin{abstract}
We give an alternative formulation of Stein's maximal inequality for generalised spherical averages in terms of derivatives of standard spherical means: if \[ k \ge 0, \qquad d \ge 2 k + 3 , \qquad \frac{d}{d - k - 1} < p < \frac{d - 1}{k} , \] and $\sigma$ is the normalised surface measure on the unit sphere $\sph$, then the maximal operator \[f \mapsto \sup_{r > 0} \, \biggl\lv r^k (\tfrac{d}{dr})^k \int_{\sph} f(\cdot + r y) \sigma(dy) \biggr\rv\] is bounded on $L^p$, with a constant that is independent of the dimension $d$.
\end{abstract}

\maketitle
\thispagestyle{empty}

%
%

\section{Introduction}

Let $\sigma$ be the normalised surface measure on the unit sphere $\sph$ in the Euclidean space $\R^d$. For a Schwartz function $f$ on $\R^d$ and $r > 0$, we define the spherical averages
\formula{
 A_r f(x) & = \int_{\sph} f(x + r y) \sigma(dy) ,
}
and if $k = 0, 1, 2, \ldots$\,, we write
\formula{
 A_r^{(k)} f(x) & = r^k (\tfrac{d}{dr})^k A_r f(x) \\
 & = \int_{\sph} r^k (\tfrac{d}{dr})^k f(x + r y) \sigma(dy) .
}
The corresponding maximal function is then defined as
\formula{
 A_*^{(k)} f(x) & = \sup \{ \lv A_r^{(k)} f(x) \rv : r \in (0, \infty) \} .
}
The classical spherical maximal function, introduced by Stein in~\cite{stein}, corresponds to $k = 0$. For $k > 0$, $A_r^{(k)} f$ can be expressed in terms of Stein's generalised spherical means, also studied in~\cite{stein}, and the maximal inequality for the latter easily leads to the following result.

\begin{theorem}[a variant of Stein's maximal inequality]
\label{thm:stein}
Suppose that
\formula[eq:range]{
 k & = 0, 1, 2, \ldots, \qquad d \ge 2 k + 3 , \qquad \frac{d}{d - k - 1} < p < \frac{d - 2}{k} ,
}
with the convention that $(d - 2) / k = \infty$ if\/ $k = 0$. Then $A_{\smash{*}}^{(k)}$ satisfies the strong type inequality
\formula[eq:stein]{
 \lV A_*^{(k)} f \rV_p & \le C_{k, p} \lV f \rV_p ,
}
with a constant $C_{k, p}$ that is independent of the dimension $d$. Consequently, $A_{\smash{r}}^{(k)}$ and $A_{\smash{*}}^{(k)}$ extend continuously to all of\/ $L^p(\R^d)$.
\end{theorem}

Although the maximal inequality for spherical means $A_r f$ (Theorem~1 in~\cite{stein}) is well known, the estimate for Stein's generalised spherical means (Theorem~2 therein) is much less familiar. The aim of this note is to present this latter result in a more accessible form and to draw broader attention to it.

Our interest in derivatives of the spherical means $A_r f$ and the associated maximal functions $A_{\smash{*}}^{(k)} f$ stems from the study of maximal inequalities for truncated Riesz transforms in~\cite{kkw}. In the companion paper~\cite{kk}, Theorem~\ref{thm:stein} leads to a more direct proof of the dimension-free estimate of Kucharski, Wróbel and Zienkiewicz from~\cite{kwz}, which improved special cases covered in~\cite{kw,lmz}, and a dimension-dependent inequality from~\cite{mopv,mov,mv}.

First order derivatives $A_{\smash{r}}^{(1)} = r \tfrac{d}{dr} A_r$ appear naturally in the study of the spherical maximal operator $A_{\smash{*}}^{(0)}$, also on spaces more general than $\R^d$; see, for example, Section~6.1 in~\cite{bhrt} by Bagchi, Hait, Roncal and Thangavelu on spherical means on the Heisenberg group. Other results partially related to the estimates of $A_{\smash{r}}^{(1)}$ include bounds for the $p$-variation of $r \mapsto A_r f$ by Jones, Seeger and Wright in~\cite{jsw}, and Beltran, Oberlin, Roncal, Seeger and Stovall in~\cite{borss}.

In~\cite{bourgain}, Bourgain proved the maximal estimate for spherical means $A_r f$ when $d = 2$ and $p > 2$. For a simplified argument and further discussion, we refer to~\cite{mss} by Mockenhaupt, Seeger and Sogge.

The left endpoint $p = d / (d - k - 1)$ in~\eqref{eq:range} is best possible, as can be seen by considering any nonnegative and nonzero $f$. On the other hand, the right endpoint is closely related to the regularity of solutions of the wave equation and Sogge's local smoothing conjecture. Miao, Yang and Zheng in~\cite{myz} observed that Bourgain and Demeter's decoupling theorem from~\cite{bd} extends the range of the maximal inequality for Stein's generalised spherical means. When $d \ge 2 k + 3$, this improves the right endpoint in Theorem~\ref{thm:stein} to $p = (d - 1) / k$, which is shown to be optimal by considering $f$ equal to mollifications of $\sigma$; a detailed discussion and further results are given by Liu, Shen, Song and Yan in~\cite{lssy}.

The above extensions of Stein's estimates lead to the following more general variant of Theorem~\ref{thm:stein}. This result is best possible when either $d \ge 2 k + 3$, or $k = 0$ and $d = 2$. The only range of parameters for which the validity of the maximal inequality~\eqref{eq:stein} remains unknown is thus $k \ge 1$, $d = 2 k + 2$, and $2 < p < 2 + \tfrac{1}{k}$.

\begin{theorem}[extension due to Bourgain, Demeter, Miao, Yang, Zheng]
\label{thm:bourgain:demeter}
The maximal estimate~\eqref{eq:stein} holds under the more general condition:
\formula[eq:range:bourgain:demeter]{
 k & = 1, 2, \ldots, \qquad d \ge 2 k + 3 , \qquad \frac{d}{d - k - 1} < p < \frac{d - 1}{k} ;
}
or
\formula[eq:range:bourgain]{
 k & = 0, \qquad d \ge 2 , \qquad \frac{d}{d - 1} < p .
}
\end{theorem}

The constants in Stein's original result, as well as its extensions discussed above, depend on the dimension. However, Stein observed in~\cite{stein:essay} that maximal estimates for certain operators improve with the dimension. While~\cite{stein:essay}, as well as the detailed proof given by Stein and Strömberg in~\cite{ss}, focused on the Hardy--Littlewood maximal function, the inequality for the spherical maximal function $A_{\smash{*}}^{(0)} f$ with constants independent from the dimension follows by the same argument, and the case of $A_{\smash{*}}^{(k)} f$ for an arbitrary $k \ge 0$ is no different.

Many subsequent works have focused on maximal inequalities for averages over dilations of a given convex body; we refer to the survey by Deleaval, Guédon and Maurey~\cite{dgm}. A variety of results also exist in other geometrical settings. We do not attempt a full literature review here, and we only mention the recent contribution by Ganguly and Ghosh~\cite{gg} for the Heisenberg group, and the survey by Bourgain, Mirek, Stein and Wróbel~\cite{bmsw} for integer lattices.

%
%

\section{Details}

Below by $C$ we denote a generic, typically large, constant. If $C$ depends on parameters, we always list them explicitly in the subscript. The value of $C$ can change even within a single equation or inequality, so, for example, $C + C = C$.


\subsection{Stein's generalised spherical averages}

In~\cite{stein}, Stein introduced a holomorphic family of generalised spherical means $M_{\smash{r}}^{(\alpha)}$ of order $\alpha \in \C$. When $f \in \schwartz(\R^d)$ (that is, $f$ is in the Schwartz class of smooth, rapidly decaying functions), $r > 0$ and $\re \alpha > 0$, then
\formula[eq:mar]{
 M_r^{(\alpha)} f(x) & = \frac{1}{\Gamma(\alpha)} \int_{\ball} f(x + r y) (1 - \lv y \rv^2)^{\alpha - 1} dy .
}
Note that $M_{\smash{r}}^{(\alpha)}$ is usually \emph{not} a true averaging operator: we have $M_{\smash{r}}^{(\alpha)} f(x) = \pi^{d / 2} \Gamma(\tfrac{d}{2} + \alpha) f(x)$ when $f$ is constant on $B(x, r)$, the ball of radius $r$ centred at $x$.

If we denote by $\xi x$ the usual dot product in $\R^d$, by $\fourier$ the Fourier transform on $\R^d$:
\formula{
 \fourier f(\xi) & = \int_{\R^d} e^{-i \xi x} f(x) dx ,
}
and by $J_\nu$ the Bessel function of the first kind, then (see p.~171 in~\cite{sw})
\formula[eq:fmar]{
 \fourier M_r^{(\alpha)} f(\xi) & = 2^{d / 2 + \alpha - 1} \pi^{d / 2} \lv r \xi \rv^{1 - d / 2 - \alpha} J_{d / 2 + \alpha - 1}(\lv r \xi \rv) \fourier f(\xi) .
}
Here, again, $f \in \schwartz(\R^d)$, $r > 0$ and $\re \alpha > 0$. However, the expression $t^{-\nu} J_\nu(t)$ is an entire function of $\nu \in \C$ for every $t > 0$, and a bounded function of $t > 0$ for every $\nu \in \C$ (see \nist{10.2}{2} in~\cite{nist}). Therefore, equation~\eqref{eq:fmar} defines a smooth function $M_{\smash{r}}^{(\alpha)} f \in L^2(\R^d)$ for every $\alpha \in \C$.


\subsection{Connection with derivatives of spherical means}

By the differentiation rule for the Bessel function (\nist{10.6}{6} in~\cite{nist}):
\formula{
 \frac{d}{d t} \bigl[t^\nu J_\nu(t)\bigr] & = t^\nu J_{\nu - 1}(t) ,
}
we have, for $f \in \schwartz(\R^d)$, $r > 0$ and $\alpha \in \C$,
\formula{
 \frac{d}{d r} \bigl[ r^{d + 2 \alpha - 2} \fourier M_r^{(\alpha)} f(\xi) \bigr] & = \frac{2^{d / 2 + \alpha - 1} \pi^{d / 2}}{\lv \xi \rv^{d + 2 \alpha - 2}} \, \frac{d}{d r} \biggl[ \lv r \xi \rv^{d / 2 + \alpha - 1} J_{d / 2 + \alpha - 1}(\lv r \xi \rv) \fourier f(\xi) \biggr] \\
 & = \frac{2^{d / 2 + \alpha - 1} \pi^{d / 2}}{\lv \xi \rv^{d + 2 \alpha - 3}} \, \lv r \xi \rv^{d / 2 + \alpha - 1} J_{d / 2 + \alpha - 2}(\lv r \xi \rv) \fourier f(\xi) \\
 & = 2 r^{d + 2 \alpha - 3} \fourier M_r^{(\alpha - 1)} f(\xi) .
}
Applying the inverse Fourier transform, we get, by the dominated convergence theorem,
\formula[eq:derivative]{
 M_r^{(\alpha - 1)} f(x) & = \frac{1}{2 r^{d + 2 \alpha - 3}} \, \frac{d}{d r} \bigl[ r^{d + 2 \alpha - 2} M_r^{(\alpha)} f(x) \bigr] .
}
Observe that, by~\eqref{eq:mar} and integration in spherical coordinates, we have
\formula{
 r^d M_r^{(1)} f(x) & = r^d \int_{\ball} f(x + r y) dy = \int_{B(x, r)} f(z) dz = \frac{2 \pi^{d / 2}}{\Gamma(\tfrac{d}{2})} \int_0^r t^{d - 1} A_t f(x) dt .
}
Combining this with~\eqref{eq:derivative} for $\alpha = 1$, we find that
\formula[eq:sphere]{
 M_r^{(0)} f(x) & = \frac{1}{2 r^{d - 1}} \, \frac{d}{d r} \bigl[ r^d M_r^{(1)} f(x) \bigr] = \frac{\pi^{d / 2}}{\Gamma(\tfrac{d}{2})} \, A_r f(x)
}
(which also follows from~\eqref{eq:fmar} and the expression for the Fourier transform of $\sigma$, see p.~154 in~\cite{sw}). Iterating~\eqref{eq:derivative} for $\alpha = 0, -1, \ldots, -k + 1$, and combining the result with~\eqref{eq:sphere}, we find that
\formula{
 M_r^{(-k)} f(x) & = \frac{1}{2^k r^{d - 2 k - 2}} \, \biggl(\frac{1}{r} \, \frac{d}{d r}\biggr)^k \bigl[ r^{d - 2} M_r^{(0)} f(x) \bigr] \\
 & = \frac{\pi^{d / 2}}{2^k \Gamma(\tfrac{d}{2}) r^{d - 2 k - 2}} \, \biggl(\frac{1}{r} \, \frac{d}{d r}\biggr)^k \bigl[ r^{d - 2} A_r f(x) \bigr] .
}
If we expand the derivatives on the right-hand side and use $A_r^{(j)} f(x) = r^j (\tfrac{d}{d r})^j A_r f(x)$, we get
\formula{
 M_r^{(-k)} f(x) & = \frac{\pi^{d / 2}}{2^k \Gamma(\tfrac{d}{2})} \biggl(A_r^{(k)} f(x) + \sum_{j = 0}^{k - 1} C_{d, k, j} A_r^{(j)} f(x) \biggr) .
}
Solving this triangular system of linear equations for $r^k A_{\smash{r}}^{(k)} f(x)$ leads to
\formula[eq:sphere:k]{
 A_r^{(k)} f(x) & = \frac{2^k \Gamma(\tfrac{d}{2})}{\pi^{d / 2}} \biggl( M_r^{(-k)} f(x) + \sum_{j = 0}^{k - 1} C_{d, k, j} M_r^{(-j)} f(x) \biggr) .
}


\subsection{Maximal inequality for Stein's generalised spherical averages}

Theorem~2 in~\cite{stein} states that for every $f \in \schwartz(\R^d)$ and $\alpha \in \C$, the maximal function
\formula{
 M_*^{(\alpha)} f(x) & = \sup\{\lv M_r^{(\alpha)} f(x) : r \in (0, \infty)\}
}
satisfies the strong type inequality on $L^p(\R^d)$:
\formula[eq:generalised]{
 \lV M_*^{(\alpha)} f \rV_p & \le C_{d, \alpha, p} \lV f \rV_p
}
whenever
\formula[eq:generalised:range:stein]{
 d & \ge 3 , \qquad \frac{d}{d + \alpha - 1} < p < \left\{\begin{aligned}
  & \frac{d - 2}{-\alpha} && \text{if } \alpha < 0 , \\
  & \infty && \text{if } \alpha \ge 0 .
 \end{aligned}\right.
}
In~\cite{bourgain} Bourgain proved~\eqref{eq:generalised} for $d = 2$ and $p > 2$. Using the decoupling theorem of Bourgain and Demeter from~\cite{bd}, Miao, Yang and Zheng in~\cite{myz} extended the range further to
\formula[eq:generalised:range]{
 d & \ge 2 , \qquad \frac{d}{d + \alpha - 1} < p < \left\{\begin{aligned}
  & \frac{2 d - 6}{-4 \alpha + 1 - d} && \text{if } \alpha < -\dfrac{(d - 1)^2}{2 d + 2} , \\
  & \frac{d - 1}{-\alpha} && \text{if } -\dfrac{(d - 1)^2}{2 d + 2} \le \alpha < 0 , \\
  & \infty && \text{if $\alpha \ge 0$.}
 \end{aligned}\right.
}
It is a challenging open problem to determine the optimal range of $p$ for which~\eqref{eq:generalised} holds.


\subsection{Maximal inequality for derivatives of spherical averages}

The relation~\eqref{eq:sphere:k} between the derivatives $A_{\smash{r}}^{(k)}$ of spherical averages and Stein's generalised spherical averages $M_{\smash{r}}^{(\alpha)}$, and the maximal inequality~\eqref{eq:generalised} for the latter, applied with $\alpha = 0, -1, -2, \ldots, -k$, immediately lead to the maximal inequality~\eqref{eq:stein} given in Theorem~\ref{thm:stein}, but with a constant $C_{d, k, p}$ that depends also on the dimension. The range of admissible parameters $d, k, p$ given in~\eqref{eq:range} follows directly from condition~\eqref{eq:generalised:range:stein} with $\alpha = -k$, as the condition then also holds for $\alpha = 0, -1, -2, \ldots, -k$.

We now argue that the extension described in Theorem~\ref{thm:bourgain:demeter} follows if we use~\eqref{eq:generalised:range} instead of~\eqref{eq:generalised:range:stein}. If~\eqref{eq:range:bourgain} holds, that is, if $k = 0$, $d = 2$ and $p > 2$, we clearly have~\eqref{eq:generalised:range} with $\alpha = -k = 0$, as desired. Suppose now that~\eqref{eq:range:bourgain:demeter} is satisfied, that is,
\formula{
 k & \ge 1 , \qquad d \ge 2 k + 3 , \qquad \frac{d}{d - k - 1} < p < \frac{d - 1}{k} .
}
We claim that condition~\eqref{eq:generalised:range} holds with $\alpha = 0, -1, -2, \ldots, -k$.

Fix $j = 0, 1, \ldots, k$. Clearly,
\formula{
 p & > \frac{d}{d - k - 1} \ge \frac{d}{d - j - 1} = \frac{d}{d + \alpha - 1} ,
}
which is the lower bound of~\eqref{eq:generalised:range}. Since $d \ge 2 k + 3$ and $(d - 3) (d + 1) \le (d - 1)^2$, we have
\formula{
 -\alpha & = j \le k \le \frac{d - 3}{2} \le \frac{(d - 1)^2}{2 d + 2} ,
}
so that the upper bound of~\eqref{eq:generalised:range} follows from
\formula{
 p & < \frac{d - 1}{k} \le \frac{d - 1}{j} = \frac{d - 1}{-\alpha} .
}
This shows our claim, and hence completes the proof of Theorem~\ref{thm:bourgain:demeter}, but with a constant $C_{d, k, p}$ in the maximal inequality~\eqref{eq:stein} that depends on the dimension.


\subsection{Independence of the dimension}
\label{sec:independence}

In order to prove that the constant in the maximal inequality~\eqref{eq:stein} in Theorems~\ref{thm:stein} and~\ref{thm:bourgain:demeter} can be chosen independently of the dimension $d$, one adapts in a straightforward way Stein's argument from~\cite{stein:essay}, described thoroughly in~\cite{ss}. The key observation is that $A_r f(x)$ can be evaluated by averaging lower dimensional spherical means of $f$ over all hyperplanes which pass through $x$. For completeness, we provide the details.

We now make the dimension explicit in the notation and write $\sph_d$, $\sigma_d$, $A_{\smash{\kern.5pt d, r}}^{\kern-.5pt(k)}$ and $A_{\smash{\kern.5pt d, *}}^{\kern-.5pt(k)}$ instead of $\sph$, $\sigma$, $A_{\smash{r}}^{(k)}$ and $A_{\smash{*}}^{(k)}$. We denote by $O_d$ the group of orthogonal $d \times d$ matrices, and by $\Sigma_d$ the normalised Haar measure on $O_d$. As in equation~(4.1) in~\cite{ss}, we have
\formula{
 \int_{\sph_d} f(y) \sigma_d(dy) & = \int_{O_d} \biggl(\int_{\sph_{d - 1} \times \{0\}} f(L y) (\sigma_{d - 1} \times \delta_0)(dy) \biggr) \Sigma_d(dL) ,
}
where $\sigma_{d - 1} \times \delta_0$ is the product measure concentrated on $\sph_{d - 1} \times \{0\}$. Applying the above identity to $f(x + r y)$ instead of $f(y)$ leads to
\formula{
 A_{d, r} f(x) & = \int_{\sph_d} f(x + r y) \sigma_d(dy) \\
 & = \int_{O_d} \biggl(\int_{\sph_{d - 1} \times \{0\}} f(x + r L y) (\sigma_{d - 1} \times \delta_0)(dy) \biggr) \Sigma_d(dL) \\
 & = \int_{O_d} \biggl(\int_{\sph_{d - 1} \times \{0\}} f \circ L(L^{-1} x + r y) (\sigma_{d - 1} \times \delta_0)(dy) \biggr) \Sigma_d(dL) \\
 & = \int_{O_d} (A_{d - 1, r} \otimes \id)(f \circ L)(L^{-1} x) \Sigma_d(dL) ,
}
where $(A_{d - 1, r} \otimes \id) f(x)$ is the operator $A_{d - 1, r}$ applied to the restriction $f(\cdot, x_d)$ of $f$ to the hyperplane $\R^{d - 1} \times \{x_d\}$ that contains $x$, evaluated at $(x_1, \ldots, x_{d - 1})$. By the dominated convergence theorem, we have
\formula{
 A_{d, r}^{(k)} f(x) & = \int_{O_d} (A_{d - 1, r}^{(k)} \otimes \id)(f \circ L)(L^{-1} x) \Sigma_d(dL) ,
}
and hence
\formula[eq:rotations]{
 \lv A_{d, *}^{(k)} f(x) \rv & \le \int_{O_d} (A_{d - 1, *}^{(k)} \otimes \id)(f \circ L)(L^{-1} x) \Sigma_d(dL) .
}
If the maximal inequality on $L^p(\R^{d - 1})$ for $A_{\smash{\kern.5pt d - 1, r}}^{\kern-.5pt(k)}$ holds with a constant $c_{d - 1, k, p}$, then, by convexity of the $L^p(\R^d)$ norm, tensorisation and invariance of the Lebesgue measure under $L$, we have
\formula*[eq:rotations:norm]{
 \lV A_{d, *}^{(k)} f \rV_p & \le \int_{O_d} \lV (A_{d - 1, *}^{(k)} \otimes \id)(f \circ L)\lV_p \Sigma_d(dL) \\
 & \le c_{d - 1, k, p} \int_{O_d} \lV f \circ L \lV_p \Sigma_d(dL) = c_{d - 1, k, p} \lV f \lV_p ,
}
and so the maximal inequality on $L^p(\R^d)$ for $A_{\smash{\kern.5pt d, r}}^{\kern-.5pt(k)}$ holds with the same constant $c_{d - 1, k, p}$.


\subsection{Fractional orders}

Iterating~\eqref{eq:derivative}, we find that if $\alpha = -k + \beta$, where $k = 0, 1, 2, \ldots$ and $\beta \in (0, 1)$, then
\formula{
 M_r^{(\alpha)} f(x) & = \frac{1}{2^k r^{d - 2 k + 2 \beta - 2}} \, \biggl(\frac{1}{r} \, \frac{d}{d r}\biggr)^k \bigl[ r^{d + 2 \beta - 2} M_r^{(\beta)} f(x) \bigr] .
}
Using~\eqref{eq:mar}, integrating in spherical coordinates, and using the dominated convergence theorem, we obtain
\formula{
 M_r^{(\alpha)} f(x) & = \frac{1}{2^k r^{d - 2 k + 2 \beta - 2}} \, \biggl(\frac{1}{r} \, \frac{d}{d r}\biggr)^k \biggl[ \frac{r^{d + 2 \beta - 2}}{\Gamma(\beta)} \int_{\ball} f(x + r y) (1 - \lv y \rv^2)^{\beta - 1} dy \biggr] \\
 & = \frac{1}{2^k r^{d - 2 k + 2 \beta - 2}} \, \biggl(\frac{1}{r} \, \frac{d}{d r}\biggr)^k \biggl[ \frac{2 \pi^{d / 2} r^{d + 2 \beta - 2}}{\Gamma(\beta) \Gamma(\tfrac{d}{2})} \int_0^1 t^{d - 1} (1 - t^2)^{\beta - 1} A_{r t} f(x) dt \biggr] \\
 & = \frac{\pi^{d / 2}}{2^{k - 1} \Gamma(\beta) \Gamma(\tfrac{d}{2})} \int_0^1 t^{d - 1} (1 - t^2)^{\beta - 1} \frac{1}{r^{d - 2 k + 2 \beta - 2}} \biggl(\frac{1}{r} \, \frac{d}{d r}\biggr)^k \bigl[ r^{d + 2 \beta - 2} A_{r t} f(x) \bigr] dt .
}
If we expand the derivatives on the right-hand side, we identify $M_{\smash{r}}^{(\alpha)} f(x)$ with a linear combination of fractional derivatives of Erdélyi--Kober type of standard spherical means $A_r f(x)$:
\formula{
 M_r^{(\alpha)} f(x) & = \sum_{j = 0}^k C_{d, \alpha, j} r^{2 j} \int_0^1 t^{d + j - 1} (1 - t^2)^{\beta - 1} A_{r t}^{(j)} f(x) dt .
}
Note that the argument involving averaging over all orthogonal transformations used in the previous section does not apply directly to~$M_r^{(\alpha)}$.


\subsection{Vector valued functions}

Suppose that $f$ takes values in some Banach space $X$ of sequences; thus, $f = (f_1, f_2, \ldots)$, where $f_1, f_2, \ldots$ are Borel functions on $\R^d$. Consider the maximal inequality of Fefferman--Stein type:
\formula[eq:fefferman:stein]{
 \biggl( \int_{\R^d} \lV A_*^{(k)} f(x) \rV_X^p dx \biggr)^{1 / p} & \le C_{d, k, p, X} \biggl( \int_{\R^d} \lV f(x) \rV_X^p dx \biggr)^{1 / p} ,
}
where we denoted $A_*^{(k)} f = (A_*^{(k)} f_1, A_*^{(k)} f_2, \ldots)$. We claim that if~\eqref{eq:fefferman:stein} holds in some dimension $d$, then it also holds in every higher dimension, with the same constant.

To prove this claim, we follow the argument from Section~\ref{sec:independence}. We apply~\eqref{eq:rotations} to $f_1, f_2, \ldots$ and use the convexity of the norm on $X$:
\formula{
 \lV A_{d, *}^{(k)} f(x) \rV_X & \le \int_{O_d} \bigl\lV (A_{d - 1, *}^{(k)} \otimes \id)(f \circ L)(L^{-1} x) \bigr\rV_X \Sigma_d(dL) .
}
Hence, denoting
\formula{
 \lV f \rV_{X, p} & = \biggl( \int_{\R^d} \lV f(x) \rV_X^p dx \biggr)^{1 / p}
}
and repeating the calculation in~\eqref{eq:rotations:norm}, we obtain
\formula{
 \lV A_{d, *}^{(k)} f \rV_{X, p} & \le \int_{O_d} \lV (A_{d - 1, *}^{(k)} \otimes \id)(f \circ L)\lV_{X, p} \Sigma_d(dL) \\
 & \le c_{d - 1, k, p, X} \int_{O_d} \lV f \circ L \lV_{X, p} \Sigma_d(dL) = c_{d - 1, k, p, X} \lV f \lV_{X, p} .
}
It follows that $c_{d, k, p, X} \le c_{d - 1, k, p, X}$, as claimed.

The above observation was suggested to the author by Błażej Wróbel. The case $k = 0$ and $X = \ell^q$ corresponds to the Fefferman--Stein inequality for spherical averages, proved under the assumption that $d \ge 3$ and $p, q \in (\tfrac{d}{d - 1}, d)$ in~\cite{dk}. Extensions to $k > 0$ are given in the companion paper~\cite{kk}.


\subsection{Remark}

The proof of the maximal inequality~\eqref{eq:generalised} for Stein's generalised spherical averages $M_{\smash{r}}^{(\alpha)}$ in~\cite{stein} can be easily adapted so that it leads to a constant $C_{d/2 + \alpha, p}$ that depends on the dimension $d$ and the order $\alpha$ only through the value of $\tfrac{d}{2} + \alpha$. As this fact does not seem to have important implications, we only sketch the necessary modification of the proof of Theorem~2 in~\cite{stein}. Instead of a compactly supported function $\ph$, we choose $\ph(x) = \pi^{-d / 2} e^{-\lv x \rv^2}$. Then, we apply Stein's maximal theorem for symmetric Markovian semigroups (with a constant independent of the dimension thanks to Rota's martingale trick and the Burkholder--Davies--Gundy inequality; see p.~106 in Section~IV.4 in~\cite{stein:topics}, or~\cite{trevisan}) instead of the estimate for the Hardy--Littlewood maximal function.

%
%

\section*{}

\subsection*{Acknowledgements}

I thank Rodrigo Bañuelos, Maciej Kucharski, Paweł Plewa, Bartosz Trojan, and Błażej Wróbel for inspiring discussions about the subject of the paper.

%
%

%
%

\end{document}